\documentclass{amsart}
\usepackage{comment}

\newtheorem{thm}{Theorem}[section]
\newtheorem{prop}[thm]{Proposition}
\newtheorem{cor}[thm]{Corollary}
\newtheorem{lem}[thm]{Lemma}

\newcommand{\R}{\mathbb{R}}
\newcommand{\N}{\mathbb{N}}
\newcommand{\Z}{\mathbb{Z}}
\newcommand{\boB}{\mathcal{B}}
\newcommand{\boL}{\mathcal{L}}
\newcommand{\dd}{\, {\rm d}}
\newcommand{\tq}{\, :\, }

\newcommand{\norm}[1]{\left\| #1 \right\|}
\newcommand{\moins}{\backslash}
\newcommand{\vf}{\varphi}
\newcommand{\au}{\underline{a}}

\renewcommand{\hat}{\widehat}

\begin{document}

\title[Regularity of coboundaries]{
Regularity of coboundaries
for non uniformly expanding Markov maps
}
\author{S\'ebastien Gou\"ezel}
\date{January 25, 2005}
\address{
D\'epartement de math\'ematiques et applications\\
\'Ecole Normale Sup\'erieure\\
45 rue d'Ulm, 75005 Paris, France.}
\email{sebastien.gouezel@ens.fr}

\begin{abstract}
We prove that solutions $u$ of the equation $f=u-u\circ T$ are
automatically H\"older continuous when $f$ is H\"older continuous, and
$T$ is non uniformly expanding and Markov. This result applies in particular to
Young towers and to intermittent maps.
\end{abstract}

\maketitle

\section{Results}

Let $(X,m)$ be a probability space and $T:X\to X$ an ergodic measure
preserving transformation. Let also $G$ be a locally compact abelian
group, endowed with an invariant metric that we denote by $|x-y|$.
It is often important to know whether a function $f:X
\to G$ is a \emph{measurable coboundary}, i.e., there exists a
measurable function $u:X \to G$ such that
  \begin{equation}
  \label{def_coboundary}
  f=u-u\circ T
  \end{equation}
almost everywhere. For $G=\R$, 
this condition is indeed often the only obstruction
to have a non-degenerate central limit theorem for the Birkhoff sums
of $f$ (see e.g.\ \cite{leonov}, \cite{guivarch-hardy}, \cite{liverani:CLT}). 
For $G=S^1$, it is relevant to prove local limit theorems
(see
\cite{aaronson_denker} and \cite{aaronson_denker_sarig_zweimuller}
when $f$ is locally constant, 
in the Markov and non-Markov case).

When $T$ is uniformly hyperbolic and $f$ is H\"older continuous, the
Liv{\v{s}}ic regularity theorem (\cite{livsic}) states that $u$ 
must have a H\"older continuous version, for which \eqref{def_coboundary}
holds everywhere. In particular, if there exists a point
$x$ such that $T^n(x)=x$ and $\sum_{k=0}^{n-1} f(T^k x) \not=0$,
then $f$ is not a measurable coboundary. Hence, it is possible to
prove in practice that a function is not a coboundary (see also
\cite{pollicott_yuri} and \cite{nicol_scott}).

In this note, we extend the aforementioned result of
Liv{\v{s}}ic to non-uniformly expanding Markov dynamical systems,
without any additional assumption on the functions $f$ or $u$. The result
will first be given in the abstract setting of \emph{Gibbs-Markov
maps} (see \cite{aaronson:book}). Applications to Young towers,
intermittent maps in dimension $1$ and positive recurrent Markov
shifts will also be described. 

The proof
is quite flexible since it is completely elementary and does not use
spectral theory. Hence, the same kind of arguments may be used in
other settings.

\subsection{Results for Gibbs-Markov maps}

In this paragraph, we will work in the setting of Gibbs-Markov maps,
defined in \cite[Section
4.7]{aaronson:book}.

Let us recall briefly the definitions. Let $(X,d,\boB,m)$ be a bounded 
metric
space endowed with its Borel $\sigma$-algebra and a probability
measure. A non singular map $T:X\to X$ is \emph{Gibbs-Markov} if there
exists a partition $\alpha$ of $X$ (modulo $0$) by sets of positive
measure, such that
\begin{enumerate}
\item For all $a \in \alpha$, $T(a)$ is a union (modulo $0$) of
elements of $\alpha$ and $T:a \to T(a)$ is invertible.
\item There exists a finite
subset $\{a_1,\dots,a_n\}$ of $\alpha$ with the following property:
for any $a\in \alpha$, there exist $i,j\in\{1,\dots,n\}$ such that
$a \subset T(a_i)$ and $a_j \subset T(a)$ (modulo $0$).
\item Expansion: there exists $\lambda>1$ such that $\forall a\in
\alpha$, for almost all $x,y\in a$, $d(Tx,Ty) \ge \lambda d(x,y)$.
\item Distortion: for $a\in \alpha$, let $g$ be the inverse of the
jacobian of $T$ on $a$, i.e.,
$g(x)=\frac{ \dd m_{|a}}{\dd (m\circ T_{|a})} (x)$ for $x\in a$. 
Then there exists $C$
such that, for all $a\in \alpha$, for almost all $x,y\in a$,
  $
  \left| 1-\frac{g(x)}{g(y)} \right| \le C d(Tx,Ty)$.
\end{enumerate}

Property (2), also known as the BIP 
(big images and preimages) 
property, is apparently stronger 
than the usual \emph{big image property} $\inf_{a\in \alpha}
m(Ta)>0$. However, when (4) is satisfied and $T$ is probability 
preserving, these two properties
are equivalent by \cite{sarig:BIP}.

Usually, Gibbs-Markov maps are endowed with a distance given by
$d(x,y)=\tau^{s(x,y)}$ where $\tau\in (0,1)$ and $s(x,y)$ is the
separation time of $x$ and $y$. We have chosen here to use a general
distance since it will be more convenient in the applications: our
main result will say that a function is Lipschitz continuous with
respect to $d$, which means that having more freedom to choose the
distance will give more precise results. In particular, when the
Gibbs-Markov map is obtained by coding another dynamical systems, it is
natural to use the distance induced by the original distance (see
Sections \ref{subsection_Young_towers} and
\ref{subsection_intermittent} for illustrations of this phenomenon).

For $a_0,\ldots,a_{n-1}\in \alpha$, let
$[a_0,\ldots,a_{n-1}]=\bigcap_0^{n-1}T^{-i}(a_i)$. It is a
\emph{cylinder}
 of length $n$.
For $f:X \to G$ and $Z \subset X$, set
  \begin{multline*}
  D f(Z)=\inf\{C >0 \tq \exists \Omega\subset Z \text{ with
}m(Z\moins\Omega)=0 \text{ such that }\\
  \forall x,y\in\Omega, |f(x)-f(y)| \le C
  d(x,y)\}.
  \end{multline*}

The main result of this note is the following theorem:
\begin{thm}
\label{thm_regul_cohom}
Let $(X,T,m,\alpha)$ be a probability preserving Gibbs-Markov map.
Let $f:X\to G$ satisfy $\sum_{a\in\alpha}m(a) Df(a)
<+\infty$. Let $u:X\to G$ be a measurable function such that
$f=u-u\circ T$ almost everywhere.

Then $\sup_{a_* \in \alpha_*} D
u(a_*)<\infty$, where $\alpha_*$ is the partition generated by the
images of the elements of $\alpha$. Moreover, the function $u$ is
essentially bounded.
\end{thm}
Remarks:
\begin{enumerate}
\item
Since $T$ is Markov, $\alpha_*$ is coarser than $\alpha$.  In
particular, $\sup_{a\in \alpha}Du(a)<\infty$, i.e., $u$ has a version
which is uniformly
Lipschitz on each element of the partition $\alpha$.
\item
The map $T$ is also Gibbs-Markov for the distance $d(x,y)^\gamma$ when
$\gamma \in (0,1]$. Hence, Theorem \ref{thm_regul_cohom} implies a
similar statement for H\"older functions.
\item
The proof will in fact show that there exists a constant $C$ depending
only on $T$ such that
  $\sup_{a_*\in \alpha_*} Du(a_*) \le C \sum_{a\in \alpha} m(a)Df(a)$.
In particular, when $f$ is constant on each element of $\alpha$, we
get $Du(a_*)=0$, i.e., $u$ is essentially constant on the elements of
$\alpha_*$. When $G=S^1$, we get a completely different proof of
\cite[Theorem 3.1]{aaronson_denker}.
\item
The proof would be easier under the stronger assumption
  \begin{equation*}
  \sup_{a\in \alpha}Df(a)<\infty.
  \end{equation*}
However, this assumption is too strong, since it is not compatible
with the induction process which will enable us to extend Theorem
\ref{thm_regul_cohom} to non uniformly expanding settings.
\end{enumerate}

In this paper, $\N=\{ n\in \Z, n\geq 0\}$ and $\N^*=\N \moins\{0\}$.

\subsection{Application to Young towers}
\label{subsection_Young_towers}
Let $(X,d,m)$ be a probability space endowed with a bounded metric
$d$. A map $T:X\to X$ is a \emph{Young tower} (\cite{lsyoung:recurrence})
if there exist integers
$R_l\in \N^*$ and a partition $\{\Delta_{k,l}\}_{l\in \N, k\in
\{0,\dots, R_l-1\}}$ of $X$ such that
\begin{enumerate}
\item
For all $l$ and $k<R_l-1$,
$T$ is a measurable isomorphism between $\Delta_{k,l}$
and $\Delta_{k+1,l}$, preserving $m$.
\item For all $l$,
$T$ is a measurable isomorphism between $\Delta_{R_l-1,l}$ and
$\Delta_0:= \bigcup_m \Delta_{0,m}$.
\item There exists $\lambda>1$ such that, for all $l$,
for all $x,y\in
\Delta_{0,l}$, $d(T^{R_l}x,T^{R_l}y)\ge \lambda d(x,y)$.
\item There exists $C>0$ such that, for all $l$ and $k<R_l$, for all
$x,y\in \Delta_{k,l}$, $d(x,y)\le C d(T^{R_l-k}x,T^{R_l-k}y)$.
\item For $x\in \Delta_{R_l-1,l}$, let $g(x)$ be the inverse of the
distortion of $T$ at $x$, i.e.,
$g(x)=\frac{\dd m_{|\Delta_{R_l-1,l}}}{\dd \left(m\circ
T_{|\Delta_{R_l-1,l}}\right)}(x)$. There exists $C>0$ such that, for all
$l$, for all $x,y\in \Delta_{R_l-1,l}$,
 $ 
  \left|1-\frac{g(x)}{g(y)}\right| \le C d(Tx,Ty)$.
\end{enumerate}
The third and fifth conditions mean that the returns to the basis are
expanding and have a controlled distortion. Hence, Young towers are a
good model for many non uniformly expanding maps: the map has good
properties, but after some waiting time which can be arbitrarily long.

\begin{thm}
\label{thm_regul_cobord_young}
Let $(X,T,m,d)$ be a Young tower, and let $f:X\to G$
satisfy
  \begin{equation*}
  \sum m(\Delta_{k,l})Df(\Delta_{k,l})<\infty.
  \end{equation*}
If $u:X\to G$ is such that $f=u-u\circ T$ almost everywhere,
then the function $u$ has a version which is Lipschitz on $\Delta_0$, i.e.,
there exists $C>0$ such that, for almost all $x,y\in \Delta_0$,
$|u(x)-u(y)| \le Cd(x,y)$.
\end{thm}

This result applies in particular when the function $f$ is Lipschitz.
\begin{proof}
By \cite{lsyoung:recurrence}, we can assume without loss of generality
that $m$ is invariant.

Let $Y=\Delta_0$ with the partition $\alpha=\{\Delta_{0,l}\}$,
$\vf:Y\to \N^*$ the first return time to $Y$ (i.e., on
$\Delta_{0,l}$, $\vf=R_l$), and $T_Y=T^\vf$ the map induced by $T$
on $Y$. Define also a distance $d'$ on $\Delta_{0,l}\in \alpha$ by
$d'(x,y)= d(T^{R_l} x,T^{R_l}y)$. If $x$ and $y$ are in two different
elements of the partition $\alpha$, set also $d'(x,y)=\lambda
\sup_{X\times X} d$. Then
$(Y,T_Y,m_{|Y}/m(Y),d')$ is a Gibbs-Markov map for the partition
$\alpha$. Moreover, $T_Y$ preserves the measure $m_{|Y}/m(Y)$ and the
partition $\alpha_*$ is the trivial partition.

Let  $f:X\to G$  satisfy  $\sum
m(\Delta_{k,l})Df(\Delta_{k,l})<\infty$, and assume that $f=u-u\circ T$.
Define a function $f_Y$ on $Y$ by
$f_Y(x)=\sum_{k=0}^{\vf(x)-1} f(T^k x)$.
On $\Delta_{0,l}$,
  \begin{align*}
  |f_Y(x)-f_Y(y)|&\le \sum_{k=0}^{R_l-1} |f(T^k x)-f(T^k y)|
  \le \sum_{k=0}^{R_l-1} d(T^k x, T^k y) Df(\Delta_{k,l})
  \\&
  \le Cd(T^{R_l}x,T^{R_l}y) \sum_{k=0}^{R_l-1} Df(\Delta_{k,l})
  = C d'(x,y)\sum_{k=0}^{R_l-1} Df(\Delta_{k,l}).
  \end{align*}
Hence, $\sum_{a\in \alpha} m(a)Df_Y(a) \le C \sum
m(\Delta_{k,l})Df(\Delta_{k,l}) < \infty$.
Moreover, $f_Y=u-u\circ T_Y$.

Theorem \ref{thm_regul_cohom} applies and proves that $u$ is almost
everywhere Lipschitz
on each element of $\alpha_*$, for the distance $d'$. 
In particular, on any element
$\Delta_{0,l}$ of $\alpha$, we get
$|u(x)-u(y)|\le E d'(x,y)$.

Take finally $x',y'\in \Delta_0$. They have preimages
$x,y$ under $T^{R_l}$ in $\Delta_{0,l}$. As $f_Y(x)=u(x)-u(x')$ and
$f_Y(y)=u(y)-u(y')$, we get
  \begin{align*}
  |u(x')-u(y')| &\le |f_Y(x)-f_Y(y)|+|u(x)-u(y)|
  \le C d'(x,y) + E d'(x,y)
  \\&
  =C' d(x',y').
  \qedhere
  \end{align*}
\end{proof}

\subsection{Applications to intermittent maps}
\label{subsection_intermittent}
For
$\alpha\in (0,1)$, let $T$ be the map from $[0,1]$ to itself given by
   \begin{equation*}
   T(x)=\left\{ \begin{array}{cl}
   x(1+2^\alpha x^\alpha) &\text{if }0\le x\le 1/2;
   \\
   2x-1 &\text{if }1/2<x\le 1.
   \end{array}\right.
   \end{equation*}
This map has been studied by \cite{
liverani_saussol_vaienti}. It is nonuniformly expanding since the
fixed point $0$ satisfies $T'(0)=1$, and admits an absolutely
continuous invariant probability measure.
\begin{prop}
\label{prop_LSV}
Let $f:[0,1]\to G$ be H\"older with exponent $\gamma>0$ on the
intervals $[0,1/2]$ and $(1/2,1]$. If $u:[0,1]\to G$ is measurable
and satisfies $f=u-u\circ T$ Lebesgue almost everywhere, then
there exists a function $\tilde{u}$, equal to $u$ almost everywhere,
H\"older with exponent $\gamma$, and such
that $f=\tilde{u}-\tilde{u}\circ T$ everywhere.
\end{prop}
\begin{proof}
Let $Y=(1/2,1]$, $\vf$ the first return time from $Y$ to itself and
$T_Y:Y\to Y$ the induced map. Then $T_Y$ is Gibbs-Markov for the
partition $B_n=\{ y\in Y \tq \vf(y)=n\}$, by \cite{
liverani_saussol_vaienti}. Hence, the arguments in the proof of Theorem
\ref{thm_regul_cobord_young} apply and prove that $u$ is a.e.\ H\"older on
$Y$. As $T: (1/2,1] \to (0,1]$ is Lipschitz and has Lipschitz inverse,
the coboundary equation implies that $u$ is a.e.\ H\"older on $(0,1]$, i.e.,
there
exists a set $V$ of full measure and a constant $C$ such that, for all
$x,y\in V$, $|u(x)-u(y)| \le C |x-y|^\gamma$.

The function $u$ is uniformly continuous on $V$, whence it can be
extended to a continuous -- and even H\"older -- function $\tilde{u}$
on $[0,1]$. On $V\cap T^{-1}(V)$, which is dense, we have
$f(x)=\tilde{u}(x)-\tilde{u}(Tx)$. Since both members of this equality are
continuous on the intervals  $[0,1/2]$ et $(1/2,1]$, this
equality holds in fact everywhere.
\end{proof}

In particular, if $f$ is a measurable coboundary, it satisfies
$\sum_0^{n-1}f(T^k x)=\tilde{u}(T^n
x)-\tilde{u}(x)=0$ at any point $x$ such that $T^n(x)=x$.
\begin{cor}
If $f:[0,1]\to \R$ is H\"older continuous on $[0,1/2]$ and $(1/2,1]$
and satisfies $f(0)\not=0$, then
$f$ is not a measurable coboundary.
\end{cor}

This solves a conjecture stated in
\cite{haydn_vaienti:entropy_fluctuation}: in this article, the
authors need to know that $f=\log |T'| -\int \log|T'|$ is not a
coboundary to get a nonzero variance in the central limit theorem.
As $f$ is $\alpha$-H\"older on $[0,1/2]$ and $(1/2,1]$, and
$f(0)=-\int \log|T'| <0$, the corollary applies and proves that it
is indeed never the case.

Using Theorem \ref{thm_regul_cohom} with $G=S^1$, we can get in the
same way a stronger result:
\begin{cor}
The function $f(x)=\log|T'|-\int \log|T'|$ can not be written as
$f=u-u\circ T+\lambda q+\mu$ almost everywhere, where $u:[0,1]\to \R$
is measurable, $q:X\to \Z$ and $\lambda,\mu \in \R$.
\end{cor}
The proof is the same, using the behavior at the fixed points $0$
and $1$ to get a contradiction. This is a strong
aperiodicity result on the function $f$. By \cite[Theorem
1.2]{gouezel:local}, it implies that $f$ satisfies a local limit
theorem when $\alpha<1/2$.

\subsection{Application to positive recurrent Markov shifts}

Let $T:X\to X$ be a positive recurrent Markov shift with H\"older 
potential, as defined in
\cite{sarig:shift}, preserving the probability measure $m$. 
The map $T$ satisfies the same assumptions as a
Gibbs-Markov map, except the BIP property. We assume also that the
distance $d$ is given by $d(x,y)=\tau^{s(x,y)}$ where $\tau\in (0,1)$
and $s(x,y)$ is the separation time of $x$ and $y$. Such maps have in
general more complicated combinatorics than Young towers, but they
enjoy uniform expansion (since $d(Tx,Ty)=\tau^{-1} d(x,y)$ for all
$x,y$ in the same element of $\alpha$) while Young towers are
expanding only after many iterates.

\begin{thm}
Let $f:X\to G$ satisfy $\sum_{a\in \alpha} m(a) Df(a)<\infty$. Let
$u:X\to G$ be a measurable function such that $f=u-u\circ T$ almost
everywhere. Then, for all $a\in \alpha$, $Du(a)<\infty$. Moreover, if
$T$ is transitive,
$\sum_{a\in \alpha} m(a) Du(a)<\infty$.
\end{thm}
\begin{proof}
For $a\in \alpha$, let $T_a$ be the map induced by $T$ on $[a]$. It is 
Gibbs-Markov. Using Theorem \ref{thm_regul_cohom}, we show as in the
proof of Theorem \ref{thm_regul_cobord_young} that $Du(a)<\infty$. If
$T$ is transitive, the proof of Lemma \ref{somme_u_finie}
applies and gives $\sum m(a)Du(a)<\infty$. 
\end{proof}

\section{Proof of Theorem \ref{thm_regul_cohom}}

A Gibbs-Markov map is \emph{transitive} if, for all $a,b\in \alpha$,
there exists $n$ such that
$b\subset T^n(a) \mod
0$. When $T$ preserves a probability measure, there exists a finite
decomposition
$\alpha=\alpha_1\cup \ldots\cup \alpha_n$ such that the image of an
element of $\alpha_i$ is contained in
$X_i=\bigcup_{a\in\alpha_i}a$, and such that
$T$ is a transitive Gibbs-Markov map on $X_i$ (\cite{aaronson:book}).
To prove the theorem, it is sufficient to prove it on each $X_i$. We
can therefore assume that $T$ is transitive. 

The main step of the proof is the following lemma:
\begin{lem}
There exists $\alpha_1\in \alpha$ such that $Du(\alpha_1)<\infty$.
\end{lem}
\begin{proof}
Let $\Phi(x)=Df(a)$ when $x\in a$. This function is integrable by
assumption. In particular, there exists a set $X_1$ of full measure
such that the Birkhoff sums $S_n \Phi(x)=\sum_{k=0}^{n-1}\Phi(T^k x)$
satisfy $S_n \Phi(x)=
O(n)$ when $x\in X_1$.

There exists $X_2$ of full measure such that, if $x\in X_2$, all its
iterates satisfy: for almost all $y$ in the same element of partition
$a$ as $T^n x$,  $|f(y)-f(T^n x)|
\le Df(a) d(y,T^n x)$.

The martingale convergence theorem implies that almost every point is
a measurable continuity point of $u$:
there exists $X_3$ of full measure such that, if $x\in X_3$ and
$a_0,a_1,\ldots$ denotes the sequence of elements of $\alpha$
containing respectively $x,Tx,\ldots$, then, for all $\varepsilon>0$,
  \begin{equation*}
  \frac{ m\{y \in [a_0,\ldots,a_{n-1}] \tq
  |u(y)-u(x)|>\varepsilon\}}
  {m[a_0,\ldots,a_{n-1}]}
  \to 0.
  \end{equation*}

As $T$ is Gibbs-Markov, all its iterates have a bounded distortion
(\cite[Proposition 4.3.1]{aaronson:book}). Hence, there exists $B>0$
such that, for any measurable set $Z$ and for any cylinder of length $k$,
  \begin{equation}
  \label{bounded_dist}
  B^{-1} \frac{m(T(a_{k-1}) \cap Z)}{m(Ta_{k-1})} \le
  \frac{m( [a_0,\ldots,a_{k-1}]\cap T^{-k}
  Z)}{m[a_0,\ldots,a_{k-1}]}
  \le B \frac{m(T(a_{k-1}) \cap Z)}{m(Ta_{k-1})}.
  \end{equation}
Since $T$ has the big image property, this implies that there exists
$B'>0$ such that
  \begin{equation}
  \label{dist}
  \frac{m( [a_0,\ldots,a_{k-1}]\cap T^{-k}
  Z)}{m[a_0,\ldots,a_{k-1}]}
  \le B'm(Z).
  \end{equation}
Let $\lambda>1$ be the expansion factor of $T$ and
let $K>0$ be large enough so that
  \begin{equation}
  \label{definit_K}
  K \log \lambda >3.
  \end{equation}
Let $\alpha_1,\ldots,\alpha_N$ be a finite number of elements of
$\alpha$ such that $m(X\moins \bigcup \alpha_i) \le \varepsilon_0$
where $\varepsilon_0$ satisfies $K \log(1-B'\varepsilon_0) \ge
-1/2$.
Write
  \begin{equation*}
  Z_n=\{ x \tq \forall n^3 \le k < n^3+\lfloor K \log n \rfloor, T^k(x) \in
  \alpha_1\cup \ldots \cup \alpha_N\}.
  \end{equation*}
Let finally $X_4$ be the set of points belonging to infinitely many
$Z_n$.
\begin{lem}
The set $X_4$ has nonzero measure.
\end{lem}
\begin{proof}
Write $A=\alpha_1\cup \ldots \cup \alpha_N$.
Let us first bound $m(Z_n)$ from below. For any cylinder
$[a_0,\ldots,a_{k-1}]$, we apply \eqref{dist} to $X\moins A$, of
measure at most $\varepsilon_0$, and we get
  \begin{equation*}
  m([a_0,\ldots,a_{k-1}] \cap T^{-k}A)
  \ge (1-B' \varepsilon_0)
  m[a_0,\ldots,a_{k-1}].
  \end{equation*}
Summing these inequalities for $a_{k-1}=\alpha_1,\ldots, \alpha_N$ yields
  \begin{equation*}
  m([a_0,\ldots,a_{k-2}]\cap T^{-k+1}A \cap T^{-k}A)
  \ge (1-B' \varepsilon_0) m([a_0,\ldots,a_{k-2}] \cap  T^{-k+1}A).
  \end{equation*}
This last term is larger than $(1-B' \varepsilon_0)^2
m[a_0,\ldots,a_{k-2}]$, again by \eqref{dist}.
We get in this way by induction
  \begin{equation*}
  m\bigl([a_0,\ldots,a_{l}]\cap T^{-l-1}A\cap
  \cdots \cap T^{-k}A\bigr) \ge (1-B'\varepsilon_0)^{k-l}
  m[a_0,\ldots,a_{l}].
  \end{equation*}
In particular, for $l=-1$ and $k=\lfloor K\log n \rfloor -1$,
we get using the invariance of $m$ that
  \begin{equation*}
  m(Z_n) \ge (1-B'\varepsilon_0)^{K \log n}=n^{K
  \log(1-B'\varepsilon_0)} \ge \frac{1}{\sqrt{n}}.
  \end{equation*}

Hence, $\sum m(Z_n)=\infty$. We will use a version of the
Borel-Cantelli Lemma to conclude. Since the sets $Z_n$ are not
independent, we will use the following version of this lemma, due to
Lamperti (\cite[Proposition 6.26.3]{borel_cantelli_dependant}):

\emph{If $\sum m(Z_n)=\infty$ and
  \begin{equation*}
  \liminf_{n\to\infty} \frac{\sum_{j,k=1}^n m(Z_j \cap
  Z_k)}{\left(\sum_{k=1}^n m(Z_k)\right)^2} <\infty
  \end{equation*}
then the set of points belonging to infinitely many $Z_n$ has nonzero
measure.}

To estimate $m(Z_j \cap Z_k)$, we will use the transfer operator
$\hat{T}$, defined on $L^2$ as the adjoint of the composition by
$T$. It acts continuously on the space $\boL$ of functions which are
bounded and Lipschitz on any element of $\alpha$. Moreover, by
\cite[Proposition 4.7.3]{aaronson:book}, there exist $M>0$ and
$\eta<1$ such that, for any $h\in \boL$,
  \begin{equation}
  \label{doeblin_fortet}
  \bigl\|\hat{T}^p h\bigr\|_{\boL} \le M (\eta^p \norm{h}_\boL+
  \norm{h}_1).
  \end{equation}

Let $\chi$ be the characteristic function of $A$, and
$\gamma_n=\prod_{0\le k < \lfloor K \log n \rfloor } \chi \circ T^k$:
hence, $m(Z_n)=\int \gamma_n \circ T^{n^3}=\int \gamma_n$, and
$m(Z_n \cap Z_p)=\int \gamma_n \circ T^{n^3}\cdot \gamma_p \circ
T^{p^3}$. The function $\chi$ belongs to $\boL$.
For $k>j$,
  \begin{equation}
  \label{jlksdfjlmsq}
  m(Z_j \cap Z_k)=\int \gamma_j\circ T^{j^3}\cdot \gamma_k \circ
  T^{k^3}
  =\int \hat{T}^{k^3-j^3}(\gamma_j) \cdot \gamma_k
  \le \norm{ \hat{T}^{k^3-j^3}(\gamma_j)}_\boL m(Z_k).
  \end{equation}

As $\hat{T}$ acts continuously on $\boL$, the function
  \begin{equation*}
  \delta_j=\hat{T}^{\lfloor K\log j \rfloor}(\gamma_j)
  =\hat{T}( \chi \hat{T}(\chi\cdots \hat{T}(\chi))\cdots))
  \end{equation*}
satisfies $\norm{\delta_j}_\boL \le (2M)^{K \log j}$.
The inequality \eqref{doeblin_fortet} applied to
$p=k^3-j^3-\lfloor K \log j \rfloor$ and $h=\delta_j$ yields
  \begin{equation}
  \label{sqjdfkljqsdklf}
  \begin{split}
  \bigl\|\hat{T}^{k^3-j^3}\gamma_j\bigr\|_{\boL}
  &
  \le M \left(\eta^{k^3-j^3-K \log j} \norm{\delta_j}_\boL +
  \norm{\delta_j}_1 \right)
  \\&
  \le M \left( \eta^{k^3-j^3-K \log j} (2M)^{K \log j}+ m(Z_j)
  \right)
  \end{split}
  \end{equation}
since $\norm{\delta_j}_1=\int \delta_j=\int \gamma_j$, for all these
functions are nonnegative.
Hence, \eqref{jlksdfjlmsq} and \eqref{sqjdfkljqsdklf} give
  \begin{equation*}
  | m(Z_j\cap Z_k)| \le M \eta^{k^3-j^3}
  (2M/\eta)^{K\log j} + M m(Z_j)m(Z_k).
  \end{equation*}

Finally,
  \begin{align*}
  \sum_{j<k \le n} m(Z_j \cap Z_k) &\le M\sum_{j<k}m(Z_j)m(Z_k)+
  M \sum_{j=1}^\infty \eta^{-j^3}(2M/\eta)^{K\log j} \sum_{k=j+1}^\infty
  \eta^{k^3}
  \\&
  \le M\left(\sum_{k\le n} m(Z_k)\right)^2 + M \sum_{j=1}^\infty
  \eta^{-j^3} (2M/\eta)^{K \log j} \sum_{l=(j+1)^3}^\infty \eta^l.
  \end{align*}
The last sum is bounded by
  \begin{equation*}
  M\sum_{j=1}^\infty \eta^{-j^3} (2M/\eta)^{K\log j}
  \frac{\eta^{(j+1)^3}}{1-\eta} <\infty,
  \end{equation*}
which shows that the aforementioned Borel-Cantelli lemma applies.
\end{proof}

We can take $x_0\in X_1 \cap X_2 \cap X_3 \cap X_4$ since this set has
positive measure. Let
$m_k\to \infty$ be such that $x_0\in Z_{m_k}$, and $n_k=m_k^3+\lfloor K
\log m_k \rfloor -1$. Then $T^{n_k}(x_0)$ belongs to one of the sets
$\alpha_1,\ldots,\alpha_N$. In particular, one of these sets is used
infinitely many times, and taking a further subsequence we can for example
assume that $T^{n_k}(x_0) \in \alpha_1$ for all $k$. We will show that
$Du(\alpha_1)<\infty$. Denote by $a_0,a_1,\ldots$ the
elements of $\alpha$ containing respectively $x_0,T(x_0),\ldots$.
Let $[\au_n]=[a_0,\dots,a_{n-1}]$, and let $v_n: Ta_{n-1} \to
[\au_n]$ be the inverse of $T^n:[\au_n] \to Ta_{n-1}$.

Let $\varepsilon>0$. As $x_0 \in X_3$,
  \begin{equation*}
  \frac{ m\{y \in [\au_{n_k}] \tq
  |u(y)-u(x_0)|>\varepsilon\}}
  {m[\au_{n_k}]}
  \to 0.
  \end{equation*}
Taking a further subsequence of $n_k$, we can assume that
  \begin{equation*}
  \sum \frac{ m\{y \in [\au_{n_k}] \tq
  |u(y)-u(x_0)|>\varepsilon\}}
  {m[\au_{n_k}]}<\infty.
  \end{equation*}
For all $k\in \N$, the distortion control \eqref{bounded_dist}
implies that
  \begin{equation*}
  \frac{m\{ y\in T a_{n_k-1} 
  \tq
  |u( v_{n_k}y)-u(x_0)|>\varepsilon\}}{m[T a_{n_k-1}]}\\
  \asymp
  \frac{m\{y\in [\au_{n_k}] \tq
  |u(y)-u(x_0)|>\varepsilon\}}{m[\au_{n_k}]}.
  \end{equation*}
Hence, $\sum_k m\{ y\in Ta_{n_k-1} \tq 
|u(v_{n_k}y)-u(x_0)|>\varepsilon\}<+\infty$. 
Therefore, $U_\varepsilon:=\{ y\in
X\tq \exists \kappa, \forall k\ge \kappa, \text{ if }y\in Ta_{n_k-1}
\text{ then }
|u(v_{n_k}y)-u(x_0)|\le\varepsilon\}$ has full measure.

Let $y_1,y_2 \in U_\varepsilon \cap \alpha_1$. If $k$ is large enough,
the preimages $y'_1$ and $y'_2$ of $y_1$ and $y_2$ under $T^{n_k}$ in
$[\au_{n_k}]$ satisfy $|u(y'_i)-u(x_0)|\le
\varepsilon$, whence $|u(y'_1)-u(y'_2)|\le 2\varepsilon$. Then
  \begin{equation}
  \label{kjqlmfjlkqsdjf}
  \begin{split}
  |u(y_1)-u(y_2)|&=
  |u\circ T^{n_k}(y'_1)-u\circ T^{n_k}(y'_2)|
  \\&
  \le \sum_{i=0}^{n_k-1} |f\circ T^i(y'_1)-f\circ T^i(y'_2)|
  +|u(y'_1)-u(y'_2)|.
  \end{split}
  \end{equation}
Recall that $n_k=m_k^3 +\lfloor K \log m_k \rfloor -1$, and that
$\Phi$ is defined by $\Phi(x)=Df(a)$ when $x\in a$. Then
  \begin{equation}
  \label{lmjkhqfdsklj}
  \begin{split}
  \sum_{i=0}^{m_k^3-1} |f\circ T^i(y'_1)-f\circ T^i(y'_2)|
  &
  \le \sum_{i=0}^{m_k^3-1} \Phi(T^i(x_0)) d(T^i y'_1, T^i y'_2)
  \\&
  \le \sum_{i=0}^{m_k^3-1} \Phi(T^i(x_0)) \lambda^{i-n_k}
  d(T^{n_k} y'_1,T^{n_k}y'_2)
  \\&
  \le \lambda^{-K \log m_k +2}  S_{m_k^3}\Phi(x_0) d(y_1,y_2).
  \end{split}
  \end{equation}
Since $x_0 \in X_1$, there exists $C$ such that $S_n
\Phi(x_0) \le Cn$ for all $n$. As $-K \log \lambda<-3$ by
\eqref{definit_K}, we get that \eqref{lmjkhqfdsklj} tends to
$0$.

Finally, set $D=\sup Df(\alpha_j)$ for $1\le j \le N$. By definition
of $m_k$, we have  $T^i(x_0) \in \bigcup_{1\le j \le N}\alpha_j$
for all $m_k^3 \le i <n_k$, whence
  \begin{align*}
  \sum_{i=m_k^3}^{n_k-1}|f\circ T^i(y'_1)-f\circ T^i(y'_2)|
  &
  \le \sum_{i=m_k^3}^{n_k-1}D d(T^i y'_1,T^i y'_2)
  \le D \sum_{i=m_k^3}^{n_k-1} \lambda^{i-n_k}d(y_1,y_2)
  \\&
  \le \frac{D}{\lambda -1}d(y_1,y_2).
  \end{align*}

Equation \eqref{kjqlmfjlkqsdjf} then yields
  \begin{equation*}
  |u(y_1)-u(y_2)|\le o(1)+\frac{D}{\lambda -1}d(y_1,y_2) +
  2\varepsilon.
  \end{equation*}

Finally, on $\alpha_1\cap \bigcap_{\varepsilon>0} U_\varepsilon$, we have
$|u(y_1)-u(y_2)|\le \frac{D}{\lambda -1}d(y_1,y_2)$.
\end{proof}

\begin{lem}
\label{somme_u_finie}
We have $\sum_{a\in \alpha}m(a) Du(a)<\infty$.
\end{lem}
\begin{proof}
Let us show that, for any $a\in\alpha$, $Du(a)<\infty$. As $T$ is
transitive, there exists $n$ such that
$a\subset T^n(\alpha_1)$. Let $[a_0,\ldots,a_{n-1}]$ be a cylinder
included in $\alpha_1$ such that
$a\subset T(a_{n-1})$. For $y_1,y_2 \in a$, let $y'_1$ and
$y'_2$ be their preimages under $T^n$ in $[a_0,\ldots,a_{n-1}]$.
Then
  \begin{equation}
  \label{kjsqlmfdjkl}
  \begin{split}
  |u(y_1)-u(y_2)|&\le \sum_{i=0}^{n-1}|f(T^i y'_1)-f(T^i
  y'_2)|+|u(y'_1)-u(y'_2)|
  \\&
  \le \sum_{i=0}^{n-1}Df(a_i) \lambda^{i-n}d(y_1,y_2)+
   \lambda^{-n} Du(\alpha_1) d(y_1,y_2),
  \end{split}
  \end{equation}
which proves that $Du(a)<\infty$.

Let $\beta$ be a finite nonempty subset of $\alpha$.
For $a\in \alpha \moins \beta$, let us show
  \begin{equation}
  \label{kljsdfklmj}
  m(a)=\sum_{n=1}^\infty 
  \sum_{a_0\in \beta, a_1,\ldots,a_{n-1}\in \alpha\moins\beta}
  m[a_0,a_1,\ldots,a_{n-1},a].
  \end{equation}
Let $Y=\bigcup_{b\in \beta} b$. Write $A_0=a$, and
$A_{n+1}=T^{-1}(A_n) \moins Y$ and $B_{n+1}=T^{-1}(A_n) \cap Y$.
We get
  \begin{equation*}
  A_n=\bigcup_{a_0,\ldots,a_{n-1}\in \alpha\moins
  \beta}[a_0,\ldots,a_{n-1},a] \text{ and }
  B_n=\bigcup_{a_0\in \beta, a_1,\ldots,a_{n-1}\in \alpha\moins \beta}
  [a_0,\ldots,a_{n-1},a].
  \end{equation*}
Thus, we want to show that $m(a)=\sum_n m(B_n)$.
The equality
$T^{-1}(A_n)=A_{n+1}\cup B_{n+1}$ implies
$m(A_n)=m(A_{n+1})+m(B_{n+1})$. By induction, we get
$m(a)=m(B_1)+\dots+m(B_n)+m(A_n)$. It remains to prove that
$m(A_n)\to 0$. Note that $A_n \subset C_n=\{ x \tq \forall
0\le k\le n, T^k(x)\not\in Y\}$. We will show that $m(C_n) \to 0$ by
proving that $C=\bigcap C_n$ has $0$ measure.
Since the measure is invariant and $C\subset T^{-1}(C)$,
$C=T^{-1}(C) \mod 0$, whence $m(C)=0$ or $1$ by ergodicity
(\cite[Theorem 4.4.7]{aaronson:book}). The set $C$ does not intersect
$Y$, which has nonzero measure, hence $m(C)=0$. This proves
\eqref{kljsdfklmj}.

Let $[a_0,\ldots,a_{n-1},a]$ be a cylinder of nonzero measure.
By
\eqref{kjsqlmfdjkl},
  \begin{equation*}
  Du(a)\le \sum_{i=0}^{n-1} \lambda^{i-n}Df(a_i)+\lambda^{-n} Du(a_0).
  \end{equation*}
Hence, \eqref{kljsdfklmj} yields
  \begin{align*}
  m(a) Du(a) \le& \sum_{n=1}^\infty 
  \sum_{a_0\in \beta, a_1,\ldots,a_{n-1}\in\alpha
  \moins \beta} m[a_0,\ldots,a_{n-1},a]\left(\sum_{i=0}^{n-1}
  \lambda^{i-n}Df(a_i)\right)
  \\& +
  m(a) \sup_{b\in \beta} Du(b).
  \end{align*}
As $\sum m(a) \sup_{b\in \beta}Du(b)<\infty$, we will show that
$\sum m(a) Du(a)<\infty$ by showing that
  \begin{equation}
  \label{xcqsdfqsfdg}
  \sum_{n=1}^\infty \sum_{a_0\in \beta, a_1,\ldots,a_{n-1}\in
  \alpha\moins \beta} m[a_0,\ldots,a_{n-1}]\left(\sum_{i=0}^{n-1}
  \lambda^{i-n}Df(a_i)\right)
  \end{equation}
is finite. In this expression, for $a'\in \alpha\moins \beta$,
the prefactor of a term $\lambda^{-k} Df(a')$ is
  \begin{multline*}
  \sum_{n=1}^\infty
  \sum_{\substack{a_0\in \beta,a_1,\ldots,a_{n-1}\in \alpha\moins
  \beta\\
  a_{n+1},\ldots,a_{n+k-1}\in
  \alpha\moins\beta}}m[a_0,\ldots,a_{n-1},a', a_{n+1},\ldots,a_{n+k-1}]
  \\
  \le
  \sum_{n=1}^\infty 
  \sum_{a_0\in \beta,a_1,\ldots,a_{n-1}\in \alpha\moins \beta}
  m[a_0,\ldots,a_{n-1},a'].
  \end{multline*}
By \eqref{kljsdfklmj}, this last term is equal to
$m(a')$. In \eqref{xcqsdfqsfdg}, the prefactor of a term
$\lambda^{-k} Df(a')$ with $a'\in \beta$ is also at most $m(a')$. Hence,
  \begin{equation*}
  \eqref{xcqsdfqsfdg}
  \le \sum_{a'\in\alpha} \sum_{k=1}^\infty m(a') \lambda^{-k} Df(a'),
  \end{equation*}
which is finite since $\sum m(a')Df(a')<\infty$.
\end{proof}

\begin{proof}[Proof of Theorem \ref{thm_regul_cohom}]

For almost all $x$, $\sum_{Ty=x} g(y)=1$. Let us write $T^{-1}(x)=\{
x_0,x_1,\ldots\}$, and let $a_i$ be the element of $\alpha$ containing 
$x_i$.
By bounded distortion and
the big image property, there exists $C>0$ such that, for all $n$,
$g(x_n) \le C m(a_n)$. As $\sum g(x_n)=1$, this implies
 $C\sum m(a_n) \ge 1$.

Let $a_*$ be an element of $\alpha_*$. Let
$x,y\in a_*$. By definition of $\alpha_*$, their preimages
$x_0,x_1,\ldots$ and $y_0,y_1,\ldots$ belong to the same elements
$a_0,a_1,\ldots$ of $\alpha$.
Since $f=u-u\circ T$, we have for any $n$
  \begin{align*}
  |u(x)-u(y)| &\le |f(x_n)-f(y_n)|+|u(x_n)-u(y_n)|
  \le (Df(a_n)+Du(a_n)) d(x_n,y_n)
  \\&
  \le (Df(a_n)+Du(a_n))\lambda^{-1} d(x,y).
  \end{align*}
Hence,
  \begin{align*}
  |u(x)-u(y)| &\le C\sum m(a_n) |u(x)-u(y)|
  \\&
  \le C \sum m(a_n) (Df(a_n)+Du(a_n)) \lambda^{-1} d(x,y).
  \end{align*}
Finally, $Du(a_*) \le \frac{C}{\lambda} \sum_{a\in \alpha}
m(a)(Df(a)+Du(a))$, which is finite by Lemma \ref{somme_u_finie}.

To prove that $u$ is essentially bounded, we use the big preimage 
property. Let 
$a_1,\dots,a_n\in \alpha$ be such that every element of $\alpha$ is
contained in the image of some $a_i$. Let $a\in \alpha$, and let $i$ be
such that $a\subset T(a_i)$. For $x\in a$, let $x'$ be its preimage in 
$a_i$, we get
  \begin{equation*}
  |u(x)|=|u(x')-f(x')| \leq \norm{u_{|a_i}}_\infty+\norm{f_{|a_i}}_\infty.
  \end{equation*}
This last quantity is uniformly bounded.
\end{proof}

\bibliography{biblio}
\bibliographystyle{alpha}

\end{document}